\documentclass[12pt,english]{article}
\usepackage[T1]{fontenc}
\usepackage[latin9]{inputenc}
\usepackage[letterpaper]{geometry}
\usepackage{amsmath}
\usepackage{amsthm}
\usepackage{graphicx}
\usepackage{setspace}

\makeatletter
\numberwithin{equation}{section}
\numberwithin{figure}{section}

\@ifundefined{date}{}{\date{}}
\usepackage{tikz}
\usepackage{pgf}
\usepackage{pgfplots}
\usepackage{pstricks} 

\makeatother

\usepackage{babel}
\begin{document}
\title{Constructive Mathematics}
\author{Mark Mandelkern }
\maketitle
\begin{quotation}
\textbf{Abstract.} Recent advances in constructive mathematics draw
attention to the importance of proofs with numerical meaning. 
\end{quotation}
An age-old controversy in mathematics concerns the necessity and the
possibility of constructive proofs. While moribund for almost fifty
years, the controversy has been rekindled by recent advances which
demonstrate the feasibility of a fully constructive mathematics. This
nontechnical article discusses the motivating ideas behind this approach
to mathematics and the implications of constructive mathematics for
the history of mathematics.

For over a hundred years the controversy over constructivity has been
simmering silently beneath the surface of mathematics, occasionally
erupting into full-blown battle, but never reaching a settlement.
Still the conflict continues. Today, while the vast majority in mathematics
use nonconstructive methods, a small minority persist in the struggle
to bring constructive methods into general use. 

During a previous episode in this controversy, Einstein asked, \textquotedbl What
is this frog-and-mouse battle among the mathematicians?\textquotedbl{}
{[}20, p. 187{]}. One might ask the same today. To begin with a brief
answer which is elaborated below, constructive proofs are those which
ultimately reduce to finite constructions with the integers 1,2,3,\,...\,.
Such proofs are said to have \emph{numerical meaning.} In \emph{constructive}
mathematics, numerical meaning is central. In contrast, \emph{classical}
mathematics, which is dominant today, admits methods which are not
in essence finite, with the result that classical proofs often lack
numerical meaning. 

Recent advances in constructivity, which point to a final resolution
of the problem, stem from Errett Bishop's 1967 book, \emph{Foundations
of Constructive Analysis} {[}1{]}. A growing number of mathematicians
now work in accord with the general principles proposed and developed
by Bishop.

The present article discusses the basic ideas which motivate this
group of modern constructivists. The treatment deals with no technical
details, but rather with fundamental human attitudes towards mathematics,
its meaning, and its purpose. (For a thorough technical discussion,
see {[}25{]}.) The author is indebted to Y. K. Chan, Michael Goldhaber,
Keith Phillips, and the late Errett Bishop for critical readings of
early drafts of this article.\\

\noindent \textbf{Historical background}\\

\noindent The present controversy in mathematics has traces even in
the early history of mathematics:
\begin{quotation}
\begin{flushright}
\emph{Not geometry, but arithmetic alone will provide satisfactory
proof.}\\
 Archytas, ca. 375 BCE {[}13, p. 49{]}
\par\end{flushright}

\end{quotation}
Archytas refers to an ancient controversy concerning geometric and
arithmetic proofs. Certain aspects of this controversy correspond
roughly to the present controversy between classical and constructive
proofs.

Throughout the history of mathematics, both constructive and nonconstructive
tendencies are found. The Pythagoreans tried to reduce all mathematics
to numbers. Plato introduced, through his theory of forms, an idealistic
approach to philosophical problems which pervades all classical mathematics;
he taught that truth exists independently of humans, who must seek
it through dialectic. Aristotle began the systematization of logic;
his law of excluded middle (discussed below) is a major cause of nonconstructivities
when applied to the mathematical infinite. Gauss first gave his complex
numbers a geometric representation, but later considered this inadequate
and gave an arithmetic formulation. In the late nineteenth century
a violent attack on nonconstructive methods was led by Leopold Kronecker
in Berlin; his conviction was: \textquotedbl God made the integers,
all else is the work of man\textquotedbl{} {[}13, p. 988{]}. At the
turn of the century, Henri Poincare in Paris strongly advocated constructive
methods; interested in applications, he criticized the classical approach,
saying, \textquotedbl True mathematics is that which serves some
useful purpose\textquotedbl{} {[}19{]}. 

Although classical mathematics has long been dominant, it was severely
challenged during the early part of the century. The challengers,
led by L. E. J. Brouwer in Amsterdam, were critical of current practice
and called for a new beginning, a careful reconstruction of the basic
mathematical framework. The defense, content with classical methods,
was led by David Hilbert in G{\small{}ö}ttingen.

Brouwer (1881-1966) demonstrated that classical mathematics is deficient
in numerical meaning. Beginning in 1907, he devoted much of his life
to attacking classical methods whose validity he questioned, showing
that these methods did not produce mathematical objects which are
explicitly constructed and which ultimately reduce to the integers.
Brouwer held that in the \textquotedbl constructive process\,...\,lies
the only possible foundation for mathematics\textquotedbl{} {[}13,
p. 1200{]} {[}5{]}.

Hilbert (1862-1943) was a leader in the development of classical mathematics.
His solution of \textquotedbl Gordan's Problem\textquotedbl{} in
1888 had accelerated the growth of modern mathematics. His proof,
however, was nonconstructive. Paul Gordan (1837-1912) in Erlangen,
who had for twenty years tried to solve the problem, exclaimed, \textquotedbl Das
ist nicht Mathematik. Das ist Theologie!\textquotedbl{} {[}20, p.
34{]}. Leading the classical defense against the constructivists,
Hilbert angrily complained that \textquotedbl forbidding a mathematician
to make use of the law of excluded middle is like forbidding an astronomer
his telescope or a boxer the use of his fists\textquotedbl{} {[}13,
p. 1204{]}. Hermann Weyl (1885-1955) in Zurich strongly supported
Brouwer and the constructivists. He charged nonconstructive proofs
with lack of significance and value, saying that classical analysis
is \textquotedbl built on sand\textquotedbl{} {[}13, p. 1203{]}.

The Brouwer-Hilbert debate, and the subsequent work of the Brouwerian
school, was devoted more to logical, philosophical, and foundational
considerations, than to positive constructive developments. A significant
advance in the latter direction was made in 1967 by Errett Bishop
(1928-1983) in California; his book {[}1{]} succeeds in developing
a large portion of analysis in a realistic, constructive manner.\newpage{}

\begin{figure}
\includegraphics[viewport=-200bp 0bp 103bp 109bp]{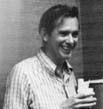}

\hspace{3.05in}{\small{}Errett Bishop}{\small\par}
\end{figure}

The central idea in modern constructive mathematics is: 
\begin{quotation}
Bishop's Thesis: \emph{All mathematics should have numerical meaning
}{[}1, p. ix{]}.
\end{quotation}
The significance of this depends upon two essential conditions. First,
classical mathematics is deficient in numerical meaning. Second, it
is in fact possible to give most mathematics numerical meaning. The
first was demonstrated by Brouwer. Although he and his followers also
made a great effort to demonstrate the second, and did constructivize
certain isolated portions of mathematics, they unfortunately introduced
unnecessary idealistic elements into much of their work. Thus Brouwer's
main contribution, of crucial significance to the development of mathematics,
was his \emph{critique of classical mathematics.} The second step,
the systematic \emph{constructive development of mathematics,} was
begun in 1967 by Bishop. A small number of workers now continue this
development, in constructive algebra {[}12{]}, {[}2l{]}, {[}24{]},
{[}26{]}; constructive analysis {[}2{]}, {[}3{]}, {[}4{]}, {[}6{]},
{[}9{]}, {[}10{]}, {[}15{]}, {[}16{]}, {[}17{]}; constructive probability
theory {[}7{]}, {[}8{]}, and constructive topology {[}22{]}. 
\begin{quotation}
Bishop's maxim: \emph{When a man proves a number to exist, he should
show how to find it.} {[}l, p.2{]} 
\end{quotation}
expresses the constructivist thesis. His book, written to demonstrate
the feasibility of \textquotedbl a straight-forward realistic approach
to mathematics,\textquotedbl{} rebuilds the basics of analysis in
a fully constructive manner and provides the framework and methods
for continuing further work \textquotedbl to hasten the inevitable
day when constructive mathematics will be the accepted norm\textquotedbl{}
{[}l. p.x{]}.\\

\noindent \textbf{Classical vs. constructive mathematics}\\

\noindent To avoid a common misunderstanding, it should be stressed
that from the constructivist position, classical mathematics does
not appear useless, but merely limited. The limitation is crucial,
but not fatal. A theorem proved by classical methods is merely incomplete;
the degree of incompleteness varies to both extremes. A classical
theorem may exhibit no numerical meaning, and there may be little
chance of extracting any numerical meaning from it. At the other extreme
is the classical theorem that is constructive as it stands, or becomes
constructive after some quite minor reformulations. Situated between
these extremes, most classical theorems are neither constructively
valid nor completely devoid of numerical meaning. Such a theorem must
be significantly modified and rephrased to show its true constructive
content, and a considerable amount of work must be done to find a
constructive proof. Often several new theorems are obtained, revealing
different constructive aspects of the classical theorem which were
hidden by its classical presentation. An example will show how a classical
theorem can lead to new constructive theorems. Consider the intermediate
value theorem which says that a continuous curve, which is somewhere
below the axis and somewhere above, must somewhere cross the axis.
This theorem is not constructively valid; there is no general finite
procedure for constructing a crossing point. Furthermore, a counterexample,
of a special type due to Brouwer, convinces us that we will never
find such a procedure. What can be done? We certainly do not want
to discard such a beautiful theorem! There are two main methods for
salvaging constructive theorems from constructively invalid classical
theorems. The first method weakens the conclusion, the second strengthens
the hypotheses. In each case we must find a constructive proof. Although
the resulting theorem does not sound as strong as the original, it
is in a deeper sense much stronger\,---\,it has numerical meaning.
Applying these methods to the continuous curve theorem, we first weaken
the conclusion. We find that we are able to construct, for any small
positive number $\varepsilon$, a point of the curve which has a distance
less than $\varepsilon$ from the axis. 

The numerical content of this constructive theorem is clear. What
was the (constructively invalid) content of the original classical
theorem? The classical theorem by itself merely makes a statement;
it says that there exists a point of the curve which lies on the axis.
In what sense, we ask, does such a point exist? Can we actually find
the point? Is there a method for constructing it? The mere statement
of the classical theorem does nothing to answer these questions. When
we examine the proof to see what is actually proved, we find the following:
Assuming that no point of the curve lies on the axis, by a certain
series of deductions we obtain a contradiction. Thus \textquotedbl existence\textquotedbl{}
in the classical theorem merely means \textquotedbl nonexistence
is contradictory.\textquotedbl{} In sharp contrast, our constructive
proof of the modified theorem contains the actual construction of
a point.

To obtain a second constructive theorem, we strengthen the hypotheses,
adding the condition that the curve is defined by a polynomial. This
extra hypothesis, while restricting the scope of the theorem, enables
us to obtain the original conclusion in full constructive force. We
can construct a point of the polynomial curve and prove that it lies
exactly on the axis. Various other constructive versions of the theorem
are also easily obtained {[}1, p. 59{]}.

The first constructivization of the intermediate value theorem, which
finds a point of the curve within $\varepsilon$ of the axis, is the
more important. This is because of its wider applicability, and because,
in general, finding solutions to problems \textquotedbl within $\varepsilon$,\textquotedbl{}
rather than \textquotedbl exactly,\textquotedbl{} is quite sufficient.
This is especially clear when one considers practical applications
of mathematics. 

What does the classical proof demonstrate? It shows that the existence
of such a continuous curve, together with a proof that no point lies
on the axis, would lead to a contradiction. This is a useful aspect
of the classical proof; we need not waste time trying to construct
such a curve. Thus, while the classical theorem may have no constructive
affirmative conclusion, it does have practical (although limited)
usefulness in directing further constructive effort. Thus, far from
being useless from a constructive viewpoint, classical mathematics
serves as an invaluable guide in building the first stages of a fully
constructive mathematics. To a large extent this has already been
done. Bishop's book contains a complete constructive development of
basic analysis: the real numbers, calculus, complex analysis, metric
spaces, measure and integration, linear spaces, and more, are constructivized.
Others have extended this work and have also begun to constructivize
algebra, topology, and probability. Yet there is still much more classical
mathematics in need of constructivization.

After the classical results in a given subject have been constructivized
as far as possible, then constructive mathematics proceeds to develop
the various constructive aspects of the subject which have been uncovered.
This enrichment of a mathematical subject is a result of distinguishing
between what is constructive and what is not, and is analogous to
the enrichment of mathematics which resulted when mathematicians began
to distinguish, for example, between infinite series which were convergent,
and those which were not.\\

\noindent \textbf{The law of excluded middle}\\

\noindent The continuous curve problem illustrates the foremost cause
of nonconstructivity and lack of numerical meaning in classical mathematics:
widespread use of the law of excluded middle, which says that any
meaningful statement is either true or false.

Great care is required with the \textquotedbl either-or\textquotedbl{}
construction, which appears in the excluded middle, and which may
be subject to varied interpretations and meanings. Suppose you are
going to lunch with a friend and you wish to know whether to take
your umbrella. If, using the law of excluded middle, your friend tells
you \textquotedbl Either it will rain this noon or it will not,\textquotedbl{}
you will find this information useless. On the other hand, if either
you are told \textquotedbl It will rain this noon,\textquotedbl{}
or you are told \textquotedbl It will not rain this noon,\textquotedbl{}
then you will have received valuable constructive information. This
example may give occasion for reflection upon the sort of information,
and the degree of certitude, that is found in mathematics, as compared
with, for example, meteorology. In fact, mathematical proofs do give
predictions. They predict that if certain numbers, with certain relationships
between them, are used in performing certain calculations, then the
resulting numbers will exhibit certain new relationships. Some statements,
such as \textquotedbl there exists a prime number between 17,000,000,000
and 17,000,000,017\textquotedbl , are in fact constructively either
true or false. This is because there is a finite procedure for deciding.
For this reason, we might be able to use an indirect proof, a \textquotedbl proof
by contradiction.\textquotedbl{} If the assumption that there is no
prime number in the specified interval leads to a contradiction, then
it would be constructively valid to conclude that there does exist
such a prime number. The essential ingredient to such an indirect
proof is the prior possession of the finite procedure which either
finds the desired prime number or proves there is none. The indirect
proof shows that the second alternative is contradictory, and thus
predicts that the finite procedure, when carried out, will lead to
the first alternative, the construction of a prime number in the specified
interval. It was for such finite situations only that Aristotle formulated
his rules of logic, especially the law of excluded middle. Nonconstructivities
arise when these rules are used indiscriminately in modern mathematics,
the science of the infinite.

Now recall our continuous curve. The statement \textquotedbl there
exists a point of the curve which lies on the axis\textquotedbl{}
admits of no finite method to determine its truth. There are infinitely
many points of the curve. Even for a single point there is no finite
method for deciding whether it lies on the axis. In contrast to the
above statement about prime numbers, here there is no prior finite
procedure for determining one of two alternatives, and thus an indirect
proof is constructively invalid.

To the constructivist, use of the law of excluded middle in infinite
situations leads only to \emph{pseudo-existence}: nonexistence is
contradictory. It is existence which stems from a construction that
is appropriate to finite humans. While the objects of constructive
mathematics are solid objects created by finite constructions, many
of the objects of classical mathematics appear as disembodied entities
born of questionable logical laws. \\

\noindent \textbf{Applications of mathematics}\\

\noindent A constructive proof, which actually constructs a point
with certain properties, has a practical advantage over a classical
proof, which merely shows that it is unthinkable that the desired
point is nonexistent.

Wandering in the Sahara, would we be content with a nonconstructive
proof of the existence of an oasis? Would our parched throats be satisfied
with a theorem which asserts that water exists somewhere in the desert,
but gives us no clue whatever as to where? Or would we prefer a constructive
drink? The classical camel tries to reassure us that there does indeed
exist an oasis. but it cannot tell us the direction or the distance.
The constructive camel, on the other hand, gives us a direction which,
while it may not lead exactly to the oasis, will lead as close as
desired. It does give us an approximate point of the compass to follow,
and an approximate distance. We might ask it to calculate a direction
so as to pass within twenty meters of the oasis. Although it could
calculate a direction  to pass within one millimeter of the exact
center of the oasis, this calculation, while still finite, would probably
take much longer, wasting precious time. In what sense does the classical
camel give the direction? The classical camel asserts that the exact
direction exists, but it would take it an infinite amount of time
to calculate it, even to find a first approximation. Unless our goatskins
happen to contain an infinite amount of water, we might find this
classical calculation a bit lengthy. Applications of mathematics are
similar to the Sahara problem. No scientist would be content to learn
that a solution to his mathematical problem exists, but that there
is no way to calculate it. Thus all experience tends to indicate that
any mathematics that is applicable must be constructive. Although
there seem to be a few applications of nonconstructive mathematics
to theoretical physics, it is likely that it will be the constructive
content of these applications which will be useful when the theory
reaches the point of experimental verification.

If it is the constructive content of mathematics which is applicable,
then since so much mathematics currently being done is nonconstructive,
why haven't users of mathematics complained? Two facts help to understand
this. First, although most mathematicians make no effort to produce
constructive results, nevertheless their results often have a very
large (and largely hidden) constructive content. It is this hidden
constructive content that is useful and applicable; a major goal of
the constructivist program is to make this hidden content explicit.
The constructive content of a classical theorem, even of its proof,
is often sufficient for applications. On the other hand, much current
mathematics is hopelessly nonconstructive; it has no numerical meaning
and no constructivizations seem possible. Such mathematics is not
being applied and will always remain inapplicable. This brings us
to the second fact, the time lag between mathematical work and its
applications. This can run to centuries, and thus isolates current
research from the test of applicability. The mathematics of previous
ages is so useful in present applications, that there is a general
belief that all current mathematical work will certainly be usable
at some time in the future. This belief may be too optimistic.\\

\noindent \textbf{Numerical meaning}\\

\noindent The suggestion that a theorem may be disputed may at first
cause some surprise. Mathematics is often presented as a prime example
of indisputable knowledge, against which other less certain forms
of knowledge are compared. When even mathematical knowledge comes
into question, it may seem that all hope is lost. Nevertheless, two
points will clarify the situation. First, there is a limited, but
crucially important, part of mathematics about which everyone agrees:
the integers. Thus, 5 + 7 = 12 is indeed a good example of indisputable
knowledge. The importance of this seemingly small part of mathematics
is that constructive mathematics attempts to build all mathematics
upon these solid integers.

Secondly, serious misunderstanding is caused by differences in interpretation
of the meanings of theorems. Superficially, the dispute sounds quite
irresolvable. The classicist has a theorem which states that a certain
point exists. The constructivist says it does not exist, and even
has a counterexample which convinces him that such a point will never
be proved to exist. A closer examination shows that the two are using
entirely different meanings of \textquotedbl exists.\textquotedbl{}
The classicist has indeed proved that the assumption that such a point
does not exist leads to a contradiction, and using the law of excluded
middle, has concluded that such a point does exist. Thus by \textquotedbl existence\textquotedbl{}
the classicist merely means \textquotedbl nonexistence is contradictory.\textquotedbl{}
On the other hand, the constructivist, when saying that a point exists,
means that a finite procedure has been given by which the point is
explicitly constructed.

Thus the dispute is resolved when we consider theorems only in conjunction
with their proofs. After all, the statement of a theorem is nothing
more than a summary of what has been demonstrated in the proof, using
concise (and often misleading) terminology. When both examine the
proof, the classicist and constructivist fully agree about what has
been proved. Has the controversy vanished into thin air? No! Rather,
we have come to the crux of the issue. The crucial question is, What
theorems should we prove? The classicist says that the theorem just
proved settles the problem, and that the constructivist is wasting
time with details. The constructivist says that the classicist has
a theorem which is splendid as far as it goes, and which points the
way to an interesting and useful constructive theorem, but which in
itself is incomplete. In the example of the continuous curve, it is
not enough to stop when the nonexistence of the point sought is shown
to be contradictory; we continue until we have constructed a point.
The construction reduces to the construction of certain integers;
it has numerical meaning. Thus the constructivist's answer to the
question, \textquotedbl What theorems should we prove?\textquotedbl ,
is given by Bishop's thesis: \textquotedbl Theorems with numerical
meaning!\textquotedbl{} It is a natural human tendency, a metaphysical
impulse, to believe that every meaningful statement must be either
true or false. This is understandable, since we are finite beings,
and usually speak only of finite matters. But in exploring the mathematical
infinite, we might heed Plutarch: \textquotedbl When talking about
infinity we are on treacherous ground and we should just try and keep
our footing\textquotedbl{} {[}18{]}. We can avoid the quicksand of
excluded middle and keep to the constructive trail. A basic metaphysical
problem is whether truth exists independently of man. The classical
approach to mathematics presumes that truth does exist in itself,
perhaps in some Platonic sphere, and we have only to find it. The
constructivist believes that mathematics belongs to man, and that
we ourselves create it, except for the integers. These integers, which
have been created for us, have already blazed the trail for us to
follow in our creation of further mathematical truth.\\

\noindent \textbf{Constructive real numbers}\\

\noindent Both geometry and arithmetic are products of human thought,
based on our concepts of space and integer. By emphasizing arithmetic
proofs we are asserting that our concept of integer is more reliable.
Even in geometry (and related fields, including analysis) we expect
more reliable results if the geometric concepts are reduced to arithmetic
concepts. We reduce the concept of a point in the plane to the concept
of real number; a point has co{\small{}ö}rdinates which are real numbers.
Then real numbers are reduced to rational numbers; a real number is
generated by a sequence of approximating rational numbers, e.g., the
finite portions of an infinite decimal. Finally, a rational number
is a ratio of integers. The classicist also reduces points and numbers
to integers in this manner. A close examination, however, shows that
the classical reduction uses the law of excluded middle, which leads
to properties of points and numbers which are constructively invalid.
A striking example of this is the trichotomy of real numbers. This
is the classical theorem that says that for any real number x, either
$x<0$, or $x=0$, or $x>0$. Constructively, it is not true; there
is no known general finite procedure which, for each real number $x$,
leads to a proof of one of the three alternatives. Worse, there is
a Brouwerian counterexample which shows that we can never expect to
find such a procedure. 

To see roughly why trichotomy fails constructively, suppose that the
real number $x=0.a_{1}a_{2}a_{3}\thinspace...$ is given in decimal
form. In any finite length of time we can calculate only finitely
many of the digits $a_{i}$. At some point we might have calculated
a million digits and found them all zero. Still, we cannot in general
predict whether all the potentially calculable digits will be zero,
or whether a nonzero digit might someday appear. Thus we cannot tell
whether $x=0$, or $x>0$. This example also indicates why the continuous
curve theorem discussed above fails constructively. We have no finite
process to decide whether the ordinate of a given point is zero or
not; we can only calculate approximations to it. Thus, while we can
tell whether a point of the curve lies near the axis, we usually cannot
tell whether it lies exactly on the axis. In the Sahara, the classical
camel might tell us to head west if $x=0$, but to head east if $x>0$.
Since there are infinitely many decimal digits, we would need to perform
an infinitely long calculation before we could take even the first
step.

The idea of the trichotomy certainly arises intuitively when we draw
a 1ine, mark a zero point on it, and look at the picture. However,
this geometric picture of the real numbers, though useful, can be
misleading. If we look at a line we may be tempted to think only of
points which we deliberately mark off on it, and for which we have
a preconceived notion about whether they lie to the left, right, or
at the zero point. Thus the trichotomy may seem evident, but this
naive view does not take into account all real numbers which have
already been constructed, or which may be constructed in the future.

The inadmissibility of trichotomy may seem a mighty blow to constructivism,
since it may seem to be such a fundamental property of the real numbers.
On the contrary, our fondness for the trichotomy arises only from
its habitual use; it is not essential for constructive analysis. In
its place we use other properties of the real numbers which are constructively
valid. For instance, given any small positive number $\varepsilon$,
it is constructively true that every real number is either less than
$\varepsilon$ or greater than zero; we have a dichotomy with a small
overlap. This constructive dichotomy is used in the construction of
a point on our continuous curve within a distance $\varepsilon$ from
the axis. \\

\noindent \textbf{The relation of constructive mathematics to the
whole of mathematics}\\

\noindent Some classical mathematicians view constructive mathematics
as a special, minor part of mathematics which, for questionable reasons,
avoids the use of certain logical principles and methods of proof.
Some followers of Brouwer maintain that constructive mathematics forms
a separate branch of mathematics, alongside of and distinct from classical
mathematics. The position of modern constructivists differs from each
of these. To them it is classical mathematics which is part of the
totality of mathematics; this totality is constructive mathematics.
The part which is classical, very large today, but very small in the
inevitable future, is that part which uses, as an extra hypothesis,
the law of excluded middle. This extra hypothesis has its limited
use, as noted above, but its general use is totally unwarranted.

Often constructive mathematics is considered part of formal logic,
philosophy, or the foundations of mathematics. Our point of view,
however, is that formal logic and philosophy make only an attempt
to lay a \textquotedbl foundation\textquotedbl{} for mathematics.
The practice of mathematics requires no foundation other than that
it be based on finite constructions which ultimately reduce to the
integers. In the recent history of constructivism, two main questions
arise. After Brouwer's critique of classical mathematics, why was
nonconstructive mathematics not immediately rejected? And now, with
the basic part of mathematics already constructivized, and with methods
for further constructive progress at hand, why do only a few use these
constructive methods? The answer to the first question is that although
Brouwer's critique of classical mathematics was clear, compelling,
widely discussed, and accepted as devastating, still no one had a
solution to the problem. Brouwer showed the lack of numerical meaning
in classical mathematics, but did not convincingly show how mathematics
might be done constructively. In fact, Brouwer himself and his followers
were convinced that it was not possible to rebuild most mathematics
constructively. They thought that most of the beautiful structures
of mathematics would be necessarily lost, and they were willing to
suffer this loss for the sake of constructibility {[}11, p. 11{]}.
It was not until later, in 1967, that Bishop showed that there need
be no loss, but rather a gain of clarity, precision, applicability,
and beauty.

The second question is much more difficult. We can only tentatively
suggest some possible reasons for the slow growth of constructivism
in mathematics since 1967. Some inertia in the university system and
the graduate curriculum may be a factor. Connected with this is the
heavy burden on students to master an extensive curriculum; this does
not leave sufficient time for exploring new ideas. Another factor
may be the chilling effect that the political troubles of this century
have had on independent thought. The slow growth of constructive mathematics
may also be related to the decay of the idea that the purpose of mathematics
is to serve the sciences. Finally, many have been diverted from a
careful investigation into the meaning of their own work by the existence
of those branches of mathematics which attempt to lay a foundation
for the rest of mathematics, validate it, and give it meaning. It
is not at all clear that this attempt will succeed. If meaning is
to be found in a piece of mathematics, that is where it will be found.
When you prove a theorem, you must yourself show where the meaning
lies; you cannot leave this task to others. 

Truth is said to reside in a deep well {[}14{]}. Reach not for the
jug of excluded middle to slake a false thirst; strive to draw what
is truly needed from the well. \\

\noindent \textbf{References}\\

\begin{singlespace}
{[}1{]} E. Bishop, \emph{Foundations of Constructive Analysis,} McGraw-Hill,
New York, 1967.

{[}2{]} E. Bishop and H. Cheng, \emph{Constructive measure theory,}
Memoirs Amer. Math. Soc., 116, 1972.

{[}3{]} D. S. Bridges, \emph{Constructive Functional Analysis,} Pitman,
London, 1979. 

{[}4{]} D. S. Bridges, A constructive look at positive linear functionals
on \emph{L(H),} Pacific J. Math., 95 (1981) 11-25.

{[}5{]} L. E. J. Brouwer, \emph{Collected Works, Vol. 1, Philosophy
and Foundations of Mathematics,} A. Heyting, Ed., North-Holland, Amsterdam,
1975.

{[}6{]} Y. K. Chan, A constructive study of measure theory, Pacific
J. Math., 41 (1972) 63-79.

{[}7{]} Y. K. Chan, A constructive approach to the theory of stochastic
processes, Trans. Amer. Math. Soc., 165 (1972) 37-44.

{[}8{]} Y. K. Chan, Notes on constructive probability theory, Ann.
Probability, 2 (1974) 51-75. 

{[}9{]} Y. K. Chan, Foundations of constructive potential theory,
Pacific J. Math, 7I (1977) 405-418. 

{[}10{]} H. Cheng, A constructive Riemann mapping theorem, Pacific
J. Math., 44 (1973) 435-454.

{[}11{]} A. Heyting, \emph{Intuitionism, 2nd ed. rev.,} North-Holland,
Amsterdam, 1966.

{[}l2{]} W. Julian, R. Mines, and F. Richman, Algebraic numbers, a
constructive development, Pacific J. Math., 74 (1978) 91-102. 

{[}13{]} M. Kline, \emph{Mathematical Thought from Ancient to Modern
Times,} Oxford University Press, New York, 1972.

{[}14{]} J. Larmor, in his 1905 introduction to the English translation
of H. Poincar{\small{}é,} \emph{Science and Hypotheses.} Reprint,
Dover, New York, 1952, p. xx. 

{[}15{]} M. Mandelkern, Connectivity of an interval, Proc. Amer. Math.
Soc., 54 (1976) 170-l72. 

{[}16{]} M. Mandelkern, Located sets on the line, Pacific J. Math.,
95 (1981) 401-409. 

{[}17{]} M. Mandelkern, \emph{Constructive continuity,} Memoirs Amer.
Math. Soc., 277, 1983. 

{[}18{]} Plutarch, Decline of the Oracles, sec. 37 in \emph{Moral
Essays,} translated by R. Warner, Penguin, 1911, p. 19.

{[}19{]} H. Poincar{\small{}é}, \emph{Science and Method,} ca. 1900.
Reprint of the English translation, ca. 1905, Dover, New York, 1952,
p.189. See also {[}13, p. 1198{]}. 

{[}20{]} C. Reid, \emph{Hilbert,} Springer, Berlin, 1970. (This book
contains much of interest for the early history of constructive mathematics.)
\-

{[}21{]} F. Richman, Constructive aspects of Noetherian rings, Proc.
Amer. Math. Soc., 44 (1974) 436-44l. 

{[}22{]} F. Richman, G. Berg, W. Julian, and R. Mines, The constructive
Jordan curve theorem, Rocky Mountain J. Math., 5 (1975) 225-236. 

{[}23{]} M. Rosenblatt, ed., \emph{Errett Bishop: Reflections on Him
and His Research,} Contemporary Mathematics. Amer. Math. Soc., 1985. 

{[}24{]} G. Stolzenberg, Constructive normalization of an algebraic
variety, Bull. Amer. Math. Soc., 74 (1968) 595-599. 

{[}25{]} G. Stolzenberg, Review; \emph{Foundations of constructive
analysis,} by E. Bishop, Bull. Amer. Math. Soc., 76 (1970) 301-323. 

{[}26{]} J. Tennenbaum, A constructive version of Hilbert's basis
theorem, dissertation, University of California at San Diego, 1973.
\\

\end{singlespace}

\noindent October 22, 1984 
\end{document}